\documentclass[11pt]{amsart}

\usepackage{amsmath}
\usepackage{amsthm}
\usepackage{graphicx}
\usepackage{color}
\usepackage{amsfonts}
\usepackage{amssymb}
\usepackage{psfrag}
\usepackage{amscd}
\usepackage{url}

\newcommand{\re}{\mathbb{R}}
\newcommand{\cpx}{\mathbb{C}}

\renewcommand{\P}{\mathbb{P}}

\newcommand{\lmd}{\lambda}

\newcommand{\eps}{\epsilon}

\newcommand{\dt}{\delta}

\def\af{\alpha}
\def\bt{\beta}

\newcommand{\sig}{\sigma}

\newcommand{\reff}[1]{(\ref{#1})}
\newcommand{\pt}{\partial}

\newcommand{\mc}[1]{\mathcal{#1}}

\newcommand{\bdes}{\begin{description}}
\newcommand{\edes}{\end{description}}

\newcommand{\bal}{\begin{align}}
\newcommand{\eal}{\end{align}}

\newcommand{\bnum}{\begin{enumerate}}
\newcommand{\enum}{\end{enumerate}}

\newcommand{\bit}{\begin{itemize}}
\newcommand{\eit}{\end{itemize}}

\newcommand{\bea}{\begin{eqnarray}}
\newcommand{\eea}{\end{eqnarray}}
\newcommand{\be}{\begin{equation}}
\newcommand{\ee}{\end{equation}}

\newcommand{\baray}{\begin{array}}
\newcommand{\earay}{\end{array}}

\newcommand{\bsry}{\begin{subarray}}
\newcommand{\esry}{\end{subarray}}

\newcommand{\bca}{\begin{cases}}
\newcommand{\eca}{\end{cases}}

\newcommand{\bcen}{\begin{center}}
\newcommand{\ecen}{\end{center}}

\newcommand{\bbm}{\begin{bmatrix}}
\newcommand{\ebm}{\end{bmatrix}}

\newcommand{\bmx}{\begin{matrix}}
\newcommand{\emx}{\end{matrix}}

\newcommand{\bpm}{\begin{pmatrix}}
\newcommand{\epm}{\end{pmatrix}}

\newcommand{\btab}{\begin{tabular}}
\newcommand{\etab}{\end{tabular}}

\theoremstyle{plain}
\newtheorem{theorem}{Theorem}[section]
\newtheorem{prop}[theorem]{Proposition}

\newtheorem{corollary}[theorem]{Corollary}

\theoremstyle{definition}

\newtheorem{exm}[theorem]{Example}

\newtheorem{remark}[theorem]{Remark}

\renewcommand{\subsection}[1]{
    \stepcounter{subsection}
    \settowidth{\hangindent}{\bf\thesubsection.~}
    \hangafter=1
    \bigskip\noindent
    {\bf\hbox{\thesubsection.~}#1}\par
    \nobreak
    \medskip
}
\usepackage[top=1.5in,bottom=1in,left=1.2in,right=1.2in]{geometry}

\begin{document}

\title[Algebraic Degree of Polynomial Optimization]
{Algebraic Degree of Polynomial Optimization}

\author{ Jiawang Nie \qquad  }
\address{Department of Mathematics, UC San Diego, 9500 Gilman Drive,
La Jolla, CA 92093, USA}
\email{njw@math.ucsd.edu}

\author{ Kristian Ranestad}
\address{Department of Mathematics, University of Oslo,
PB 1053 Blindern, 0316 Oslo, Norway}
\email{ranestad@math.uio.no}

\begin{abstract}
Consider the polynomial optimization problem
whose objective and constraints are all described by multivariate polynomials.
Under some genericity assumptions, 
we prove that the optimality conditions always hold on optimizers,
and the coordinates of optimizers are algebraic functions of the coefficients of the input polynomials.
We also give a general formula for the algebraic degree of the optimal coordinates.
The derivation of the algebraic degree is equivalent to
counting the number of all complex critical points.
As special cases, we obtain the algebraic degrees
of quadratically constrained quadratic programming (QCQP),
second order cone programming (SOCP)
and $p$-th order cone programming (pOCP),
in analogy to the algebraic degree
of semidefinite programming \cite{NRS}.
\end{abstract}

\keywords{
algebraic degree, polynomial optimization,
optimality condition,
quadratically constrained quadratic programming (QCQP), $p$-th order cone programming,
second order cone programming (SOCP), variety
}

\maketitle

\section{Introduction}

Consider optimization problem
\begin{align} \label{POP}
\left\{ \baray{rl}
\underset{x\in\re^n}{\min} & \quad f_0(x) \\
s.t. & \quad f_i(x) = 0, \, i=1,\cdots,m_e \\
& \quad f_i(x) \geq 0, \, i=m_e+1,\cdots, m
\earay \right.
\end{align}
where $f_i(x)$ are multivariate polynomial functions in $\re[x]$
(the ring of polynomials in $x=(x_1,\cdots,x_n)$ with real coefficients).
The recent interest on solving polynomial
optimization problems
\cite{Las01,NDS,Par00,ParStu}
by using semidefinite relaxations or other algebraic methods
motivates this study of the algebraic properties
of the polynomial optimization problem \reff{POP}.
A fundamental problem about \reff{POP} is
how the optimal solutions depend on the input polynomials $f_i(x)$.
When the optimality condition holds and
it has finitely many complex solutions,
the optimal solutions are algebraic functions of the coefficients of polynomials $f_i(x)$,
i.e., the coordinates of optimal solutions are roots of some univariate polynomials
whose coefficients are functions of the input data.
An interesting and important problem in optimization theory
is to study the properties of these algebraic functions, e.g.,
how big their degrees are, i.e., what is the number of complex solutions to the critical equations of \reff{POP}.
Let us begin discussions with some special cases.

The simplest case of \reff{POP} is the linear programming (LP),
i.e., all polynomials $f_i(x)$ have degree one.
In this case, the problem \reff{POP} has the form
(after removing the linear equality constraints)
\begin{align} \label{LP}
\left\{\baray{rl}
\underset{x\in\re^n}{\min} & \quad c^Tx \\
s.t. & \quad Ax \geq b
\earay \right.
\end{align}
where $c,A,b$ are matrices or vectors of appropriate dimensions.
The feasible set of \reff{LP} is now a polytope described by some linear inequalities.
As is well-known, one optimal solution $x^*$ (if it exists)
of \reff{LP} must occur at one vertex of the polytope.
So $x^*$ can be determined by the linear system
consisting of the active constraints.
When the objective $c^Tx$ is changing,
the optimal solution might move from one vertex to
another vertex.
So the optimal solution is a piecewise linear fractional
function of the input data $(c,A,b)$.
When $c,A,b$ are all rational, an optimal solution
must also be rational,
and hence its algebraic degree is one.

A more general convex optimization
which is a proper generalization of linear programming is
semidefinite programming (SDP) which has the standard form
\begin{align} \label{SDP}
\left\{\baray{rl}
\underset{x\in\re^n}{\min} & \quad c^Tx  \\\
s.t.  & \quad   A_0 + \overset{n}{\underset{i=1}{\sum}} x_i A_i \succeq 0
\earay\right.
\end{align}
where $c$ is a constant vector and the $A_i$ are constant symmetric matrices.
The inequality $X\succeq 0$ means the matrix $X$ is positive semidefinite.
Recently, Nie et al.~\cite{NRS} studied
the algebraic properties of semidefinite programming.
When $c$ and $A_i$ are generic,
the optimal solution $x^*$ of \reff{SDP} is shown \cite{NRS} to be
a piecewise algebraic function of $c$ and $A_i$.
Of course, the constraint of \reff{SDP}
can be replaced by the nonnegativity of all the principle minors
of the constraint matrix,
and hence \reff{SDP} becomes a special case of \reff{POP}.
However, the problem \reff{SDP} has very special nice properties, e.g.,
it is a convex program and the constraint matrix is linear
with respect to $x$.
Interestingly, if $c$ and $A_i$ are generic,
the degree of each piece of this algebraic function
only depends on the rank of the constraint matrix at the optimal solution.
A formula for this degree is given in \cite{NRS}.

Another optimization problem frequently used
in statistics and biology is
the Maximum Likelihood Estimation (MLE),
which has the standard form
\be \label{MLE}
\max_{x\in \Theta} \,\,
 p_1(x)^{u_1} p_2(x)^{u_2} \cdots p_n(x)^{u_n}
\ee
where $\Theta$ is an open subset of $\re^n$,
the $p_i(x)$ are polynomials such that $\sum_i p_i(x) =1$,
and the $u_i$ are given positive integers.
The optimizer $x^*$ is an algebraic function of $( u_1,\ldots,u_n)$.
This problem has recently been studied and a formula for the degree of
this algebraic function has been found (cf. \cite{CHKS, HKS}).

In this paper we consider the general optimization problem \reff{POP}, when the polynomials
$f_{0},f_{1},...,f_{m}$ define a {\it complete intersection}, i.e.,
their common set of zeros has codimension $m+1$.
We show that an optimal solution is an algebraic function of the input data.
We call the degree of this algebraic function the {\it algebraic degree}
of the polynomial optimization problem \reff{POP}.
Equivalently, the algebraic degree equals the number of complex solutions to
the critical equations of \reff{POP}, when this is finite.
Under some genericity assumptions, we give in this paper a formula for the algebraic degree of \reff{POP}

Throughout this paper, the words ``generic" and ``genericity"
are frequently used. These words are given a precise meaning in algebraic geometry.
Some property or condition holds ``generically' means
it holds in some Zariski open set
(a set described by
polynomial inequalities $\ne$).
Any statement that is proved under such a
genericity hypothesis will be valid for all data that lie in a dense, open subset
of the space of data, and hence it will hold except on a set of Lebesgue measure
zero.

The algebraic degree of polynomial optimization \reff{POP}
addresses the computational complexity at a fundamental level.
To solve \reff{POP} exactly essentially reduces to
solving some univariate polynomial equations
whose degrees are the algebraic degree of \reff{POP}.
As we can see later,
the algebraic degree might be very big.

The paper is organized as follows.
Section~2 derives a general formula for the algebraic degree,
and Section~3 gives the formulae of the algebraic degrees for special cases like
quadratically constrained quadratic programming,
second order cone programming, and $p$-th order cone programming.

\section{A general formula for algebraic degree}
\setcounter{equation}{0}

In this section, we shall derive
the formula for the algebraic degree of polynomial optimization problem \reff{POP},
when the polynomials define a complete intersection.
Suppose the polynomial $f_i(x)$ has degree $d_{i}$.
Let $x^*$ be one local or global optimal solution of \reff{POP}.

At first, we assume all the inequality constraints
are active, i.e., $m_e=m$,
and the coefficients of polynomials $f_1,f_2,\cdots, f_m$ are generic.
When $m=n$, by Bertini's Theorem \cite[\S17.16]{Ha}, the feasible set of \reff{POP}
is finite and hence the algebraic degree
is equal to the B\'{e}zout's number $d_1d_2\cdots d_m$.
So, without loss of generality, assume $m<n$.
If the variety
\[
V=\{x\in \cpx^n:\, f_1(x) = \cdots = f_m(x) =  0 \}
\]
is smooth at $x^*$, i.e., the gradient vectors
\[
 \nabla f_1(x^*), \nabla f_2(x^*), \cdots,  \nabla f_m(x^*)
\]
are linearly independent,
then the {\it Karush-Kuhn-Tucker} (KKT) condition holds at $x^*$
(Chapter~12 in \cite{NW}), i.e.,
\begin{align} \label{KKT}
\left\{ \baray{r}
 \nabla f_0(x^*) + \overset{m}{\underset{i=1}{\sum}} \lmd_i^* \nabla f_i(x^*) =0 \\
f_1(x^*) = \cdots = f_m(x^*) =  0
\earay \right.
\end{align}
where $\lmd_1^*,\cdots,\lmd_m^*$ are Lagrange multipliers
for constraints $f_1(x)=0,\cdots, f_m(x)=0$.
Thus the optimal solution $x^*$ and Lagrange multipliers
$\lmd^*=(\lmd_1^*,\cdots,\lmd_m^*)$ are
determined by the polynomial system \reff{KKT}.
The set of points $x^*$ in solutions to \reff{KKT} forms the locus of critical points of \reff{POP}.
If the system \reff{KKT} is zero-dimensional, then,
by {\it elimination theory} \cite{CLO}, the coordinates of the points $x^*$
are algebraic functions of the coefficients of the polynomials $f_i$.
Each coordinate $x_i^*$ can be determined by some univariate polynomial equation like
\[
(x_i^*)^{\dt_i} + a_1 (x_i^*)^{\dt_i-1}+ \cdots + a_{\dt_i-1} x_i^*+ a_{\dt_i} = 0
\]
where $a_j$ are rational functions of the coefficients of the $f_i$.
Interestingly, when $f_1,f_{2},...,f_{m}$ are generic,
the KKT condition always holds at any optimal solutions,
and the degrees $\dt_i$ are equal to each other.
This common degree counts the number of solutions to \reff{KKT}, i.e.,
the cardinality of the critical locus of \reff{POP} or, by definition,
the algebraic degree of the polynomial optimization \reff{POP}.
We will derive a general formula for this degree.

In what follows, we work on the complex projective spaces,
where the above question may be answered as a problem in intersection theory.
For this we need to translate the optimization problem to a relevant intersection problem.
Let $\P^n$ be the $n$-dimensional complex projective space.
A point $\tilde x \in \P^n$ is a class of
vectors $(x_0,x_1, \cdots,x_n)$ that are parallel to each other.
A variety in $\P^n$ is a set of points $\tilde x$
that satisfy a collection of homogeneous polynomial equations in $(x_0,x_1,\cdots,x_n)$.
Let $\tilde f_i(\tilde x) = x_0^{d_i}f_i(x/x_0)$ be the homogenization of $f_i(x)$.
Define $\mc{U}$ to be the projective variety in $\P^n$ as
\[
\mc{U} =\{ \tilde x\in \P^n:\,
\tilde f_1(\tilde x) = \tilde f_2(\tilde x) = \cdots = \tilde f_m(\tilde x) =  0 \}.
\]
Next, we let
\[
\tilde\nabla\tilde f_i(\tilde x)=
\bbm \frac{\partial}{\partial x_{0}}\tilde f_{i} & \cdots & \frac{\partial}{\partial x_{n}}\tilde f_{i} \ebm^T
\]
be the gradient vector, with respect to the homogeneous coordinates.
Notice that $(\frac{\partial}{\partial x_{j}}\tilde f_{i}=x_0^{d_i-1}\frac{\partial}{\partial x_{j}}f_{i}(x/x_{0})) $,
so the homogenization of $\nabla f_{i}$ coincides with the last $n$ coordinates in $\tilde\nabla\tilde f_i$.

In this homogeneous setting, the optimality condition for problem~\reff{POP} with $m=m_{e}$ is
\begin{align} \label{HOP}
\left\{ (x,\mu) \in\re^n \times \re :\,
\baray{l}
\tilde f_0(\tilde x)-\mu x_{0}^{d_{0}}=\tilde f_{1}(\tilde x)= \cdots = \tilde f_{m}(\tilde x)=0\\
\mbox{rank} \big[\tilde\nabla (\tilde f_0(\tilde x)+\mu x_{0}^{d_{0}}), \tilde\nabla (\tilde f_1(\tilde x)),...,
\tilde\nabla (\tilde f_m(\tilde x))\big]\leq m
\earay \right\}
\end{align}
where $\mu\in \re$ is the critical value.
Let $\tilde x^*\in \{x_{0}\not=0\}$ be a critical point, i.e., a solution to \reff{HOP}.
We may eliminate $\mu$ by asking that the matrix
\[ \bbm
\tilde f_0(\tilde x^*) &\tilde f_1(\tilde x^*) & \cdots & \tilde f_m(\tilde x^*)  \\
x_{0}^{d_{0}} & 0& \cdots & 0 \\
\ebm \]
have rank one, and the matrix
\[ \bbm
\frac{\pt}{\pt x_0}\tilde f_0(\tilde x^*) & \frac{\pt}{\pt x_0}\tilde f_1(\tilde x^*) & \cdots &
\frac{\pt}{\pt x_0}\tilde f_m(\tilde x^*) &(d_{0}-1) x_{0}^{d_{0}}  \\
\frac{\pt}{\pt x_1}\tilde f_0(\tilde x^*) & \frac{\pt}{\pt x_1}\tilde f_1(\tilde x^*) & \cdots &
\frac{\pt}{\pt x_1}\tilde f_m(\tilde x^*) & 0 \\
\vdots & \vdots & \vdots & \vdots & \vdots \\
\frac{\pt}{\pt x_n}\tilde f_0(\tilde x^*) & \frac{\pt}{\pt x_n}\tilde f_1(\tilde x^*) & \cdots &
\frac{\pt}{\pt x_n}\tilde f_m(\tilde x^*) & 0
\ebm \]
have rank $m+1$.
The first condition and the condition $x_{0}\not=0$ mean that our critical points
\[
\tilde x^*\in \mc{U}=\{\tilde f_1(\tilde x) = \cdots = \tilde f_m(\tilde x)=0\}
\]
while the rank of the second matrix equals $m+1$ at points where  $x_{0}\not=0$  only if the submatrix
\[ M=\bbm
\frac{\pt}{\pt x_1}\tilde f_0(\tilde x) & \frac{\pt}{\pt x_1}\tilde f_1(\tilde x) & \cdots &
\frac{\pt}{\pt x_1}\tilde f_m(\tilde x)  \\
\vdots & \vdots & \vdots & \vdots  \\
\frac{\pt}{\pt x_n}\tilde f_0(\tilde x) & \frac{\pt}{\pt x_n}\tilde f_1(\tilde x) & \cdots &
\frac{\pt}{\pt x_n}\tilde f_m(\tilde x)
\ebm \]
has rank $m$.
Therefore we define
$\mc{W}$ to be the projective variety in $\P^n$:
\[
\mc{W} =
\left\{ \tilde x\in \P^n:\,
\mbox{ all the $(m+1)\times (m+1)$ minors of $M$ vanish }  \right\},
\]
the locus of points where the rank of $[\nabla ( \tilde f_{0}),...,\nabla (\tilde f_{m})]$ is less than or equal to $m$.
Denote the class of $(1,x_1,\cdots,x_n)$ in $\P^n$ by $\tilde x$.

\begin{prop} \label{prop:kkt}
Assume $m=m_e$.
If the polynomials $f_1,\cdots,f_m$ are generic, then we have:
\bit
\item [(i)] The affine variety
$V=\{x\in \cpx^n:\, f_1(x) = \cdots = f_m(x) =  0 \}$
is smooth.
\item [(ii)] The KKT condition holds at any optimal solution $x^*$;
\item [(iii)] If $f_0$ is also generic,
the affine variety
\be \label{eq:Kdef}
K = \left\{ x\in V:\, \exists \, \lmd_1,\cdots, \lmd_m  \mbox{ such that }
\nabla f_0(x) + \overset{m}{\underset{i=1}{\sum}} \lmd_i \nabla f_i(x) =0 \right\}
\ee
defined by KKT system~\reff{KKT} is finite.
\eit
\end{prop}
\begin{proof}
(i)  When polynomials $f_1,\cdots,f_m$ are generic,
by Bertini's Theorem \cite[\S 17.16]{Ha},
the variety $\mc{U}$ has codimension $m$ and is smooth, in particular the affine subvariety
$V=\mc{U}\cap \{x_{0}\not=0\}$ is smooth.
In terms of the Jacobian matrix
\[ \bbm
\frac{\pt}{\pt x_0}\tilde f_1(\tilde x) & \frac{\pt}{\pt x_0}\tilde f_2(\tilde x) & \cdots &
\frac{\pt}{\pt x_0}\tilde f_m(\tilde x)  \\
\frac{\pt}{\pt x_1}\tilde f_1(\tilde x) & \frac{\pt}{\pt x_1}\tilde f_2(\tilde x) & \cdots &
\frac{\pt}{\pt x_1}\tilde f_m(\tilde x) \\
\vdots & \vdots & \vdots & \vdots \\
\frac{\pt}{\pt x_n}\tilde f_1(\tilde x) & \frac{\pt}{\pt x_n}\tilde f_2(\tilde x) & \cdots &
\frac{\pt}{\pt x_n}\tilde f_m(\tilde x)
\ebm, \]
its rank is full at $\tilde x$.  Furthermore, the tangent space of $V$ at $\tilde x$ is, of course,
not contained in the hyperplane ${x_{0}=0}$ at infinity,
so the column $\bbm 1 & 0 & \cdots & 0 \ebm^T$ is not in the column space of the matrix at $\tilde x$.
Therefore already the submatrix
\[ \bbm
\frac{\pt}{\pt x_1}\tilde f_1(\tilde x) & \frac{\pt}{\pt x_1}\tilde f_2(\tilde x) & \cdots &
\frac{\pt}{\pt x_1}\tilde f_m(\tilde x) \\
\vdots & \vdots & \vdots & \vdots \\
\frac{\pt}{\pt x_n}\tilde f_1(\tilde x) & \frac{\pt}{\pt x_n}\tilde f_2(\tilde x) & \cdots &
\frac{\pt}{\pt x_n}\tilde f_m(\tilde x)
\ebm \]
has full rank at $\tilde x$, i.e., the gradients
\[
 \nabla \tilde f_1(\tilde x) \,\,\,  \cdots \,\,\, \nabla \tilde f_m(\tilde  x)
\]
are linearly independent at $\tilde x\in V$.

(ii) When $x^*$ is one optimizer,  which must belong to $V$, by (i),
we know the gradients
\[
 \nabla f_1(x^*), \nabla f_2( x^*), \cdots,  \nabla f_m( x^*)
\]
are linearly independent.
Hence the KKT condition holds at $x^*$ (Chapter~12 in \cite{NW}).

(iii) We claim that the intersection $\mc{U}\cap\mc{W}$ defined above is finite.
Since our critical points $V\cap\mc{W}$ is a subset of $\mc{U}\cap\mc{W}$, (iii) would follow.
The codimension of $\mc{U}$ is $m$, and this variety is smooth,
so the matrix $M$ has by (i) rank at least $m$ at each point of $\mc{U}$.
The variety $\mc{U}\cap \{\tilde f_{0}(\tilde x) =0\}$,  is, by Bertini's Theorem, also smooth, so as above,
the matrix $M$ has full rank at points in the affine part $V\cap \{f_{0}(x) =0\}$.
On the other hand, $M$ is the Jacobi matrix for the variety
$\mc{U}\cap \{\tilde f_{0}(\tilde x) =0\}$.
This variety is again smooth and has codimension $m+1$ in the hyperplane $\{x_{0}=0\}$,
so $M$ must have full rank $m+1$ on $\mc{U}\cap \{\tilde f_{0}(\tilde x) =0\}$.
The variety $\mc{W}$ where $M$ has rank at most $m$,
therefore cannot intersect $\mc{U}\cap \{\tilde f_{0}(\tilde x) =0\}$.
But B\'{e}zout's Theorem \cite[\S 8.4]{Ful}
says that if the sum of the codimensions of two varieties in $\P^n$ does not exceed $n$, then they intersect.
In particular, any curve in $\mc{U}$ intersects the hypersurface $\{\tilde f_{0}(\tilde x)=0\}$.
Since $\mc{U}\cap \{\tilde f_{0}(\tilde x)=0\}$ has codimension $m+1$,
we deduce that $\mc{W}$ must have codimension at least $n-m$.
Furthermore, since any curve in $\mc{U}\cap\mc{W}$ would intersect $\{\tilde f_{0}(\tilde x)=0\}$,
the intersection $\mc{U}\cap\mc{W}$ must be empty or finite.
On the other hand, the variety of $n\times (m+1)$-matrices with homogeneous forms
as entries having rank no more than $m$ has codimension at most $n-m$.
So the codimension of $\mc{W}$ equals $n-m$.
Hence $\mc{U}$ and $\mc{W}$ have complementary dimensions.
Therefore the intersection $\mc{U}\cap\mc{W}$ is non-empty and (iii) follows.
\end{proof}

By Proposition~\ref{prop:kkt},
for generic $f_1,\cdots,f_m$
the optimal solutions of \reff{POP}
can be characterized by the KKT system \reff{KKT},
and for generic objective function $f_0$
the KKT variety $K$ is finite.
Geometrically, the algebraic degree of the optimization problem \reff{POP}
is, under these genericity assumption, equal to the number of
distinct complex solutions of KKT, i.e., the cardinality of the variety $K$
which we above showed to coincide with $V\cap \mc{W}$.
The variety  $\mc{U}\cap\mc{W}$ above clearly contains $K$.
On the other hand, $\mc{U}\cap\mc{W}$ is finite and does not intersect the hyperplane $\{x_{0}=0\}$
when polynomials $f_i$ are generic.
Since $\mc{U}-V=\mc{U}\cap \{x_{0}=0\}$ and the $\mc{U}\cap \mc{W} \cap \{x_{0}=0\} =\emptyset$,
we can see that the cardinality of $K$ coincides with the
cardinality, i.e., the degree of $\mc{U}\cap\mc{W}$.

For integers $(n_1,n_2,\cdots,n_k)$, define
the symmetric sum of products as follows
\be \label{eq:D}
D_r(n_1,n_2,\cdots,n_k) = \sum_{i_1+i_2+\cdots+i_k = r} n_1^{i_1}\cdots n_k^{i_k}.
\ee

\begin{theorem} \label{algdeg}
Assume $m=m_e$.
If the polynomials $f_0,f_1,\cdots,f_m$ are generic,
then the algebraic degree of \reff{POP} is
\[
d_1d_2\cdots d_m D_{n-m} (d_0-1,d_1-1,\cdots,d_m-1).
\]
Furthermore, if some $f_i$ is not generic and the system~\reff{KKT} is zero-dimensional,
then the above formula is an upper bound of the algebraic degree.
\end{theorem}
\begin{proof}
When $f_1,f_2,\cdots,f_m$ are generic, $\mc{U}$ is a smooth complete intersection of codimension $m$.
Its degree $\deg(\mc{U})=d_1d_2\cdots d_m$.
When $f_0$ is also generic,
$\mc{W}$ has codimension $n-m$ and intersects $\mc{U}$ in a finite set of points as shown above.
If the intersection $\mc{U}\cap\mc{W}$ is transverse (i.e., smooth) and hence consists of a collection of simple points,
then the degree $\deg\big(\mc{U}\cap\mc{W}\big)$ counts the number of intersection points of $\mc{U}\cap\mc{W}$,
and hence the cardinality of KKT variety $K$,
which is also the number of solutions to the KKT system~\reff{KKT} for problem~\reff{POP}.

To show that this intersection is transversal,
we consider the subvariety $X$ in $\P^n\times \P^m$ defined by the $m$ equations
$\tilde f_{1}=\tilde f_{2}=...=\tilde f_{m}=0$ and the $n$ equations
$$M\cdot(\lambda_{0},...,\lambda_{m})^T=0,$$
where the $\lambda_{i}$ are homogeneous coordinate functions in the second factor.
The image under the projection of the variety $X$ defined by these $m+n$ polynomials
into the first factor coincides with the finite set $ \mc{U} \cap \mc{W}$.
Since $M$ has rank at least $m$ at every point of $\mc{U}$,
there is a unique $\tilde\lambda = (\lambda_{0},...,\lambda_{m}) \in \P^m$
for each point $\tilde x\in \mc{U} \cap \mc{W} $
such that $(\tilde x,\tilde\lambda)$ lies in $X$.
Therefore the $X$ is a complete intersection.
It is easy to check that this complete intersection does not
have any fixed point when the coefficients of $f_{0}$ varies.
So Bertini's Theorem \cite[\S17.16]{Ha} applies to conclude that
for generic $f_{0}$ this complete intersection is transversal,
which implies that the intersection $ \mc{U} \cap \mc{W}$ in $\P^n$ is also transversal.

Since the intersection $ \mc{U} \cap \mc{W}$ is finite, i.e.,
has codimension in $\P^n$ equal to the sum of the codimensions of $\mc{U}$ and $\mc{W}$,
B\'{e}zout's Theorem (cf. \cite[\S 8.4]{Ful}, \cite[Theorem 18.3]{Ha}) applies to compute the degree
\[
\deg( \mc{U} \cap \mc{W} ) = \deg( \mc{U} ) \cdot \deg( \mc{W} ).
\]

To complete the computation, we therefore need to find $\deg(\mc{W})$.  Since the codimension of $\mc{W}$
equals the codimension of the variety defined by the $(m+1)\times (m+1)$ minors of a general $n\times (m+1)$
matrix with polynomial entries, the formula of Thom-Porteous-Giambelli \cite[\S 14.4]{Ful} applies to compute this degree:
The degree of $\mc{W}$ equals the degree of the determinantal
variety of $n\times (m+1)$ matrices of rank at most $m$,
in the space of matrices whose entries in the $i$-th column are generic forms of degree $d_{i}-1$.
These matrices may be considered as a collection of linear maps parameterized by $\P^n$.
More precisely, they define a map between vector bundles of rank $m+1$ and $n+1$ over $\P^n$
\[
\mc{M}:  {\mathcal O}_{\P^n}(-d_{0}+1)\oplus  {\mathcal O}_{\P^n}(-d_{1}+1)\oplus \cdots  \oplus{\mathcal O}_{\P^n}(-d_{m}+1)
\longrightarrow {\mathcal O}_{\P^n}^{n+1},
\]
and $\mc{W}\subset \P^n$ is the variety of points over which the map has rank $m$.
The Thom-Porteous-Giambelli formula computes the degree in terms of the topological Chern classes of these vector bundles:
The degree equals the degree of the Chern class
\[
c_{n-m}\Big(\big({\mathcal O}_{\P^n}^{n+1}\big)-\big({\mathcal O}_{\P^n}(-d_{0}+1)\oplus
{\mathcal O}_{\P^n}(-d_{1}+1)\oplus \cdots  {\mathcal O}_{\P^n}(-d_{m}+1)\big)\Big)
\]
which coincides with the coefficient of $t^{n-m}$ in
\begin{align*}
\frac{1}{(1-(d_{0}-1)t)\cdot ... \cdot (1-(d_{m}-1)t)} &= \\
(1+(d_{0}-1)t+(d_{0}-1)^2t^2+\cdots)\cdots(1 +(d_{m}-1)t &+(d_{m}-1)^2t^2+\cdots).
\end{align*}
Thus $\deg(\mc{W})$ is the complete homogeneous symmetric function of degree
$\mbox{codim}(\mc{W})$ evaluated at the column degree of $G$, which is
$D_{n-m} (d_0-1,d_1-1,\cdots,d_m-1)$.
Therefore the degree formula for the critical locus $\mc{U} \cap \mc{W}$ and hence the algebraic degree of \reff{POP} is proved.

\smallskip

When some polynomial $f_i$ is not generic, then a perturbation argument can be applied.
Let $x^*$ be one fixed optimal solution of optimization problem \reff{POP}.
Apply a generic perturbation $\Delta_\eps f_i$ to each $f_i$ so that
$(f_i+\Delta_\eps f_i)(x)$ is a generic polynomial
and the coefficients of $\Delta_\eps f_i$ tends to zero
as $\eps \to 0$. Then one optimal solution $x^*(\eps)$
of the perturbed optimization problem \reff{POP} tends to $x^*$.
By genericity of $(f_i+\Delta_\eps f_i)(x)$, we know
\[
a_0(\eps)(x_i^*(\eps))^{\dt} + a_1(\eps) (x_i^*(\eps))^{\dt-1}+ \cdots + a_{\dt-1}(\eps) x_i^*(\eps)+ a_{\dt}(\eps) = 0.
\]
Here $\dt=d_1d_2\cdots d_mD_{n-m} (d_0-1,d_1-1,\cdots,d_m-1)$
and $a_j(\eps)$ are rational functions of
the coefficients of $f_i$ and $\Delta_\eps f_i$.
Without loss of generality,
we can normalize $a_j(\eps)$ such that
\[
\max_{0\leq j\leq \dt} \,\,\, | a_j(\eps) | = 1.
\]
When $\eps \to 0$, by continuity, we can see that
$x_i^*$ is a root of some univariate polynomial
whose degree is at most $\dt$ and
coefficients are rational functions of the coefficients of
polynomials $f_0, f_1, \cdots, f_m$.
\end{proof}

\begin{remark}\label{BAS}
The genericity assumption in the theorem is used to conclude that
the critical locus $\mc{U} \cap \mc{W}$ is a smooth $0$-dimensional variety by
appealing to Bertini's Theorem \cite[\S17.16]{Ha}.
A sufficient condition for Bertini's Theorem to apply can be expressed in terms of
the sets $U_{i}$ of polynomials in
which the polynomials $f_{0},f_{1},\ldots,f_{m}$ can be freely chosen.
First assume that the generic polynomial in each $U_{i}$ is reduced, and that $U_{i}$
intersects every Zariski open set of a complex affine space $V_{i}$.
Secondly, assume that the set of common zeros of all the polynomials in $\cup_{i=0}^mV_{i}$ is empty.
Then Bertini's Theorem applies.
In fact, the polynomials $f_{i}$ for which the conclusion of Bertini's Theorem fails are contained
in a complex subvariety of $V_{i}$.

If some of the polynomials $f_{i}$ are reducible, then we may replace $f_{i}$ by the factor of least degree that contains the optimizer.
The original problem \reff{POP}, is then modified to one with a smaller algebraic degree.
This is relevant in the above context, if the generic polynomial in $U_{i}$ is reducible.

\end{remark}

\begin{exm}
Consider the following special case of problem \reff{POP}
{\small
\begin{align*}
f_0(x) & =  21 x_2^2-92 x_1  x_3^2-70 x_2^2 x_3-95 x_1^4-47 x_1  x_3^3+
51 x_2^2 x_3^2+47 x_1^5+5 x_1  x_2^4+33 x_3^5, \\
f_1(x) & =  88 x_1 +64 x_1  x_2-22 x_1  x_3-37 x_2^2+68 x_1  x_2^2 x_3-84 x_2^4+
80 x_2^3 x_3+23 x_2^2 x_3^2-20 x_2 x_3^3-7 x_3^4, \\
f_2(x) & = 31-45 x_1  x_2+24 x_1  x_3-75 x_3^2+16 x_1^3-44 x_1^2 x_3
-70 x_1  x_2^2-23 x_1  x_2 x_3-67 x_2^2 x_3-97 x_2 x_3^2.
\end{align*}
}
Here $m=m_e=2$.
By Theorem~\ref{algdeg}, the algebraic degree of the optimal solution is bounded by
\[
4\cdot 3 \cdot D_1(4,3,2) = 12 \cdot (4+3+2) = 108.
\]
Symbolic computation shows
the optimal coordinate $x_1$ is a root of the univariate polynomial of degree $108$
(whose coefficients are modulo $17$)
{\tiny
\begin{align*}
&x_1^{108} +8 x_1^{107} +7 x_1^{106} +4 x_1^{105} -x_1^{104} -x_1^{103} +2 x_1^{102} -7 x_1^{100} -7 x_1^{99} +7 x_1^{98}
+5 x_1^{95} -4 x_1^{94} -6 x_1^{93} +4 x_1^{92} -8 x_1^{91} +6 x_1^{90} \\
&+4 x_1^{89} +6 x_1^{88} +2 x_1^{87} +6 x_1^{86} -7 x_1^{85} -3 x_1^{84} +5 x_1^{83} -6 x_1^{82} -3 x_1^{81} +8 x_1^{80}
-4 x_1^{79} -x_1^{78} -2 x_1^{77} +x_1^{76} -3 x_1^{75} +6 x_1^{74} \\
&+7 x_1^{73} +4 x_1^{72}+3 x_1^{71} -4 x_1^{70} -8 x_1^{68} -x_1^{67} -x_1^{66} +2 x_1^{65} +6 x_1^{64} -4 x_1^{63} +5 x_1^{62}
+2 x_1^{61}+4 x_1^{60} -2 x_1^{59} -5 x_1^{58} +7 x_1^{57} \\
&-8 x_1^{56}+5 x_1^{55} +8 x_1^{54}-8 x_1^{53} -2 x_1^{52} -2 x_1^{51} -4 x_1^{50} -3 x_1^{49} +5 x_1^{48}-6 x_1^{46} +6 x_1^{45} -6 x_1^{44}
+5 x_1^{43} +5 x_1^{42}-5 x_1^{41} -x_1^{40}\\
&+5 x_1^{39}-4 x_1^{38}-3 x_1^{37} +5 x_1^{36}-2 x_1^{35} -x_1^{34} -6 x_1^{33} -8 x_1^{32} +6 x_1^{31} +6 x_1^{30} +8 x_1^{29} +4 x_1^{28} -8 x_1^{27}
-5 x_1^{26}-4 x_1^{25} +2 x_1^{24} \\
& -x_1^{23} +2 x_1^{22}+3 x_1^{21}+2 x_1^{20} +4 x_1^{19} +6 x_1^{18}+5 x_1^{17} -7 x_1^{16} -2 x_1^{15} -x_1^{14} -7 x_1^{13}
+5 x_1^{12} +2 x_1^{11}-8 x_1^{10} -5 x_1^{9} -5 x_1^{8} \\
&-3 x_1^{7} -2 x_1^{5} -7 x_1^{4} -2 x_1^{3} -6 x_1^{2}-3x_1-1.
\end{align*}
}
In this case the degree bound $108$ is sharp.
\end{exm}

\bigskip

Now we consider the more general case that $m>m_e$, i.e., there are
inequality constraints.
Then a similar degree formula as in Theorem~\ref{algdeg} can be obtained,
when the active set is identified.
\begin{corollary}  \label{cor:deg}
Let $x^*$ be one optimizer and $j_1,\cdots, j_k$
be the active set of inequality constraints.
If every active $f_i$ is generic, then the algebraic degree of $x^*$ is
\[
d_1\cdots d_{m_{e}} d_{j_1}\cdots d_{j_k}
D_{n-{m_{e}}-k} (d_0-1,d_1-1,\cdots,d_{m_{e}}-1,d_{j_1}-1,\cdots, d_{j_k}-1).
\]
If some $f_i$ is not generic and the system \reff{KKT} is zero-dimensional,
then the above formula is an upper bound of the degree.
\end{corollary}
\begin{proof}
Note that $x^*$ is also an optimal solution of polynomial optimization problem
\begin{align*}
\left. \baray{rl}
\underset{x\in\re^n}{\min} & \quad f_0(x) \\
s.t. & \quad f_i(x) = 0, \, i=1,\cdots, m_e \\
& \quad f_i(x) = 0, \, i=j_1,\cdots, j_k \\
\earay \right\}.
\end{align*}
Hence the conclusion follows from Theorem~\ref{algdeg}.
\end{proof}

\bigskip

\section{Some special cases}
\setcounter{equation}{0}

In this section we derive the algebraic degree of
some special polynomial optimization problems.
The simplest special case is that all the polynomials $f_i$ in \reff{POP} have degree one,
i.e., \reff{POP} becomes one linear programming of the form \reff{LP}.
If the objective $c$ is generic, precisely $n$ constraints
will be active. So the algebraic degree is $D_0(0,0,\cdots,0)=1$.
This is consistent with we have observed in Introduction.
Now let us look at other special cases.

\subsection{Unconstrained optimization}

We consider the special case that
the problem \reff{POP} has no constraints.
It becomes an unconstrained optimization.
Its optimal solutions makes the gradient of the objective vanish.
By Theorem~\ref{algdeg}, the algebraic degree is
bounded by $D_n(d_0-1)= (d_0-1)^n$, which is exactly the B\'{e}zout's number
of the gradient polynomial system
\[
\nabla f_0(x) = 0.
\]
Since $f_{0}$ can be chosen freely among all polynomials of degree $d_{0}$, Remark~\ref{BAS}
applies to show that the degree bound above is sharp.

\begin{exm}
Consider the minimization of $f_0(x)$ given by
{\small
\begin{align*}
f_0 & = x_1^4+x_2^4+x_3^4+x_4^4+x_1^3+x_2^3+x_3^3+x_4^3
-13 x_1^2-30 x_1 x_2-9 x_1 x_3+5 x_1 x_4+11 x_2^2 \\
& \quad -3 x_3 x_2-3 x_3^2-20 x_3 x_4-13 x_2 x_4-9 x_4^2
+x_1-2 x_2+12 x_3-13 x_4.
\end{align*}
}
For the above polynomial, the algebraic degree of the optimal solution is $3^4=81$.
Symbolic computation shows
the optimal coordinate $x_1$ of $x^*$ is a root of the univariate polynomial of degree $81$
(whose coefficients are modulo $17$)
{\tiny
\begin{align*}
&x_1^{81} -x_1^{80} +6 x_1^{79} -x_1^{78} +2 x_1^{77} -2 x_1^{75} +5 x_1^{74} -2 x_1^{73} +4 x_1^{72} +8 x_1^{71} +6 x_1^{70}
- x_1^{69} +2 x_1^{68} -5 x_1^{67} +7 x_1^{66} -4 x_1^{65} -3 x_1^{64} \\
&+2 x_1^{63} +8 x_1^{62} +7 x_1^{61} +5 x_1^{60}+ 4 x_1^{59} +7 x_1^{58} -2 x_1^{57}
-8 x_1^{56} -2 x_1^{55} -8 x_1^{54} +2 x_1^{53} -8 x_1^{52} -x_1^{51} +8 x_1^{50} -4 x_1^{49}
-7 x_1^{48} -x_1^{47} \\
&+5 x_1^{46} +6 x_1^{45} -3 x_1^{44} +5 x_1^{43} -4 x_1^{42} -5 x_1^{41} +x_1^{40}
-4 x_1^{39} -3 x_1^{38} +8 x_1^{37} +4 x_1^{36} +2 x_1^{35} -3 x_1^{34} -7 x_1^{33} -4 x_1^{32} -5 x_1^{31} +4 x_1^{30} \\
&-4 x_1^{29} -6 x_1^{28} -8 x_1^{27} -5 x_1^{26} -8 x_1^{25} +6 x_1^{24} +7 x_1^{23} +2 x_1^{22} +5 x_1^{21} +x_1^{20}
-4 x_1^{19} -6 x_1^{18} +4 x_1^{17} +7 x_1^{16} +5 x_1^{15} \\
&+7 x_1^{14} -3 x_1^{12} +8 x_1^{11} -x_1^{10} -5 x_1^{9} -4 x_1^{8} -6 x_1^{7}
+7 x_1^{6} +2 x_1^{5} +3 x_1^{4} +7 x_1^{3} -3 x_1^{2}+5 x_1 -6.
\end{align*}
}
The degree bound $81$ is sharp for this problem.
\end{exm}

%
%

\subsection{Quadratic constrained quadratic programming}

Consider the special case that
all the polynomials $f_0,f_1,\cdots, f_m$ are quadratic.
Then problem \reff{POP} becomes one quadratic constrained quadratic programming (QCQP)
which has the standard form
\begin{align*}
\min_{x\in\re^n} & \quad x^TA_0x+b_0^Tx+c_0 \\
s.t. & \quad  x^TA_ix+b_i^Tx+c_i  \geq 0 , \, i=1,\cdots, \ell.
\end{align*}
Here $A_i,b_i,c_i$ are matrices or vectors of appropriate dimensions.
The objective and all the constraints are all quadratic polynomials.
At one optimal solution, suppose $m\leq \ell$ constraints are active.
By Corollary~\ref{cor:deg}, the algebraic degree is bounded by
\be \label{qcqpdeg}
2^m \cdot
D_{n-m}(1,\underbrace{1,\cdots,1}_{m \mbox{ times } })
 = 2^m \cdot \sum_{i_0+i_1+i_2+\cdots+i_m = n-m} 1 = 2^m  \cdot \binom{n}{m}.
\ee
The polynomials $f_{0}, f_{1},...,f_{m}$ can be chosen freely in the space of quadratic polynomials,
so Remark~\ref{BAS} applies to show that the degree bound above is sharp.

\begin{exm}
Consider the polynomials
{\small
\begin{align*}
f_0 & = -20-27 x_1^2+89 x_1 x_2+80 x_1 x_3-45 x_1 x_4+19 x_1 x_5+42 x_1-13 x_2^2+31 x_2 x_3-79 x_2 x_4 \\
& \qquad +74 x_2 x_5-9 x_2+56 x_3^2-77 x_3 x_4-2 x_3 x_5+35 x_3+40 x_4^2-13 x_4 x_5+60 x_4+58 x_5^2-84 x_5, \\
f_1 & = 33+55 x_1^2-41 x_1 x_2+33 x_1 x_3-61 x_1 x_4+96 x_1 x_5+12 x_1+74 x_2^2-90 x_2 x_3-57 x_2 x_4 \\
& \qquad -52 x_2 x_5+51 x_2+15 x_3^2+81 x_3 x_4+87 x_3 x_5+75 x_3-10 x_4^2+58 x_4 x_5+33 x_4+83 x_5^2-23 x_5, \\
f_2 &  = 8-9 x_1^2+56 x_1 x_2-24 x_1 x_3+81 x_1 x_4+85 x_1 x_5-99 x_1-77 x_2^2-75 x_2 x_3+x_2 x_4+38 x_2 x_5\\
& \qquad +23 x_2-97 x_3^2-14 x_3 x_4-73 x_3 x_5+65 x_3+3 x_4^2-14 x_4 x_5+16 x_4+9 x_5^2-10 x_5, \\
f_3 & = 9+90 x_1^2-94 x_1 x_2-22 x_1 x_3-24 x_1 x_4+78 x_1+32 x_2^2-48 x_2 x_3-6 x_2 x_4+80 x_2 x_5-18 x_2-63 x_3^2\\
& \qquad +66 x_3 x_4-13 x_3 x_5+88 x_3-45 x_4^2-92 x_4 x_5-69 x_4-43 x_5^2+32 x_5.
\end{align*}
}
For the above polynomials, the QCQP problem is nonconvex.
We consider those local optimal solutions which make
all the three inequalities active.
By Corollary~\ref{cor:deg}, the algebraic degree of this problem is bounded by $2^m\binom{n}{m}=80$.
Symbolic computation shows
the optimal coordinate $x_1$ is a root of the univariate polynomial of degree $80$
(whose coefficients are modulo $17$)
{\tiny
\begin{align*}
&x_1^{80} -3 x_1^{79} +6 x_1^{78} +2 x_1^{77} +6 x_1^{76} -3 x_1^{75} +4 x_1^{74} -6 x_1^{73} +x_1^{72} +7 x_1^{71}
-4 x_1^{69} +4 x_1^{68} -6 x_1^{67} +x_1^{65} +5 x_1^{64} +x_1^{63} -2 x_1^{62} -6 x_1^{61} \\
&\,+8 x_1^{60} +7 x_1^{59} +x_1^{58} -7 x_1^{57} +8 x_1^{56} -5 x_1^{55} -x_1^{54} -3 x_1^{53} +x_1^{52} -5 x_1^{51} -4 x_1^{50}
+3 x_1^{49} -2 x_1^{48} -x_1^{47} -7 x_1^{46} +2 x_1^{45} +8 x_1^{44} \\
&\, +6 x_1^{43} -3 x_1^{42} +5 x_1^{41} -3 x_1^{40}+5 x_1^{39} +2 x_1^{38} +2 x_1^{37} +5 x_1^{36}
+x_1^{35} +4 x_1^{34} +4 x_1^{33} +x_1^{32} -x_1^{31} +2 x_1^{30}-x_1^{29} +4 x_1^{28}\\
&\, -2 x_1^{27} +6 x_1^{26} +6 x_1^{25} +5 x_1^{24} +3 x_1^{23} +5 x_1^{22} -x_1^{21}
+2 x_1^{20} +8 x_1^{19} -x_1^{18} +7 x_1^{17} -x_1^{15} +x_1^{14} \\
&\, +4 x_1^{13} +7 x_1^{11} -8 x_1^{10}
+3 x_1^{9} +6 x_1^{8} -x_1^{7} +8 x_1^{6} -4 x_1^{5} +8 x_1^{4} +x_1^{3}-2 x_1^{2} +7 x_1-4.
\end{align*}
}
The algebraic degree of this problem is $80$ and the bound given by formula~\reff{qcqpdeg} is sharp.
\end{exm}

\subsection{Second order cone programming}

The second order cone programming (SOCP) has the following standard form
\be  \label{SOCP}
\baray{rl}
\underset{x\in\re^n}{\min} & \quad c^Tx \\
s.t. & \quad  a_i^Tx+b_i-\|C_ix+d_i\|_2 \geq 0 , \, i=1,\cdots, \ell
\earay
\ee
where $c,a_i,b_i,C_i,d_i$ are matrices or vectors of appropriate dimensions.
Let $x^*$ be one optimizer.
Since SOCP is a convex program,
the $x^*$ must also be a global solution.
By removing the square root in the constraint, SOCP becomes the
polynomial optimization
\begin{align*}
\min_{x\in\re^n} & \quad c^Tx \\
s.t. & \quad   (a_i^Tx+b_i)^2 - (C_ix+d_i)^T(C_ix+d_i) \geq 0 , \, i=1,\cdots, \ell.
\end{align*}
Without loss of generality,
assume that the constraints with indices $1,2,\cdots,m$ are active at $x^*$.
The objective is linear but the constraints are all quadratic.
As we can see, the Hessian of the constraints has
the special form $a_ia_i^T - C_i^TC_i$.
Let $r_i$ be the number of rows of $C_i$.
When $r_i=1$, the constraint
$ a_i^Tx+b_i-\|C_ix+d_i\|_2 \geq 0 $
is equivalent to two linear constraints
\[
-(a_i^Tx + b_i) \leq  C_i x + d_i \leq  a_i^Tx + b_i.
\]
Thus, when every $r_i=1$, the problem reduces to a linear programming
and hence has algebraic degree one,
because in this situation the polynomial
$(a_i^Tx+b_i)^2 - (C_ix+d_i)^2 $ is reducible.
When $r_{i}\geq 2$ and $a_{i}, b_{i}, C_{i}, d_{i}$ are generic,
the polynomial
$(a_i^Tx+b_i)^2 - (C_ix+d_i)^T(C_ix+d_i)$ is quadratic of rank $r_i+1$ and hence irreducible.
Without loss of generality,
assume $1=r_{1}=r_{2}=...=r_{k}<r_{k+1}\leq ...\leq r_{m}$.
Then problem~\reff{SOCP} is reduced to
\begin{align*}
\min_{x\in\re^n} & \quad c^Tx \\
s.t. & \quad   a_i^Tx+b_i + \sig_i (C_ix+d_i) \geq 0 , \, i=1,\cdots, k \\
     & \quad  (a_i^Tx+b_i)^2 - (C_ix+d_i)^T(C_ix+d_i) \geq 0 , \, i=k+1,\cdots, m
\end{align*}
where scalar $\sig_i$ is chosen such that $a_i^Tx^* + b_i + \sig_i(C_i x^* + d_i)=0$.
By Corollary~\ref{cor:deg}, the algebraic degree of SOCP in this modified form is bounded by
{\small
\be \label{eq:socpdeg}
2^{m-k} \cdot
D_{n-m}(0,\cdots,0,\underbrace{1,\cdots,1}_{m-k \mbox{ times } }) =
2^{m-k} \cdot \sum_{i_{k+1}+i_{k+2}+\cdots+i_{m} = n-m} 1
= 2^{m-k} \cdot \binom{n-k-1}{m-k-1} .
\ee
}
When $k=m$, we have already seen the algebraic degree is one.

For the sharpness of degree bound~\reff{eq:socpdeg}, we apply Bertini's Theorem following Remark~\ref{BAS}.
For every $i=k+1,\cdots,m$, define the set $U_{i}$ of polynomials as
\[
U_i = \left\{
 (a_i^Tx+b_i)^2 - \sum_{1\leq j \leq r_i}  \af_{j}^2 (C_ix+d_i)_j^2 :\,
\af_1,\cdots, \af_{r_i} \in \re
\right\}.
\]
Then define affine spaces $V_{i}$ as follows:
\begin{align*}
V_i &= \left\{
(a_i^Tx+b_i)^2 - \sum_{1\leq j \leq r_i}  \bt_j (C_ix+d_i)_j^2 :\,
\bt_1,\cdots, \bt_{r_i} \in \cpx
\right\}, \, i=k+1,\cdots, m.
\end{align*}
Then every set $U_{i}$ intersects any Zariski open subset of the affine space $V_{i}$.
On the other hand the set of common zeros of the linear polynomials
$$a_i^Tx+b_i + \sig_i (C_ix+d_i),\, i=1,...,k$$
and all the polynomials in the union $\bigcup_{i=k+1}^m V_i$
is contained in the set
\be \label{eq:Zdef}
Z = \bigcap_{i=1}^k \left\{ x\in\re^n: \,a_i^Tx+b_i+\sig_i(C_ix+d_i)=0 \right\}
\bigcap_{i=k+1}^m \left\{ x\in\re^n: \,\bmx a_i^Tx+b_i=0 \\  C_ix+d_i=0 \emx \right\}.
\ee
Therefore, for generic choices
$a_i,b_i,C_i,d_i$, if $r_{k+1}+\cdots+r_m + m > n$,
the set $Z$ is empty.
Hence Remark~\ref{BAS} applies to show that,
for generic choices of
$c,a_i,b_i,C_i,d_i$, if $r_{k+1}+\cdots+r_m +m > n$,
the algebraic degree bound
$2^{m-k} \cdot \binom{n-k-1}{m-k-1}$ is sharp.

\begin{exm}
Consider SOCP defined by polynomials
{\small
\begin{align*}
f_0 & = -x_1+6 x_2+13 x_3+11 x_4+8 x_5, \\
f_1 & = (-11 x_1-18 x_2-4 x_3+2 x_4-12 x_5+7)^2-(-4 x_1-10 x_2+20 x_3-4 x_4-9 x_5+3)^2 \\
& \quad -(-5 x_1-11 x_2+8 x_3-18 x_4+11 x_5+15)^2-(21 x_1+18 x_2-12 x_3-10 x_4-8 x_5+4)^2, \\
f_2 & = (-5 x_1-5 x_2-7 x_3-6 x_4+4 x_5+41)^2-(x_1-2 x_2+10 x_3-21 x_4-11)^2 \\
& \quad -(-12 x_1+3 x_2+16 x_3+4 x_4+x_5+9)^2-(14 x_1+20 x_2-13 x_3-7 x_4+4 x_5+2)^2, \\
f_3 & = (x_1-8 x_2+11 x_3-x_5+22)^2-(2 x_1-x_2+3 x_3-x_4-25 x_5-8)^2 \\
& \quad -(-2 x_1-17 x_3+14 x_4+4 x_5-7)^2-(x_1+12 x_2+14 x_3-6 x_4-4 x_5-10)^2.
\end{align*}
}
There are no linear constraints.
For this SOCP,
all the three inequalities are active at the optimizer.
All the matrices $C_i$ has three rows.
By formula~\reff{eq:socpdeg}, the algebraic degree of this problem is bounded by
$2^3\binom{5-1}{3-1}=48$.
Symbolic computation shows that
the optimal coordinate $x_1$
is a root of the univariate polynomial of degree $48$
(whose coefficients are modulo $17$)
{\tiny
\begin{align*}
& x_1^{48}-2x_1^{47}-3x_1^{46}+3x_1^{45}+4x_1^{44}+5x_1^{43}-6x_1^{42}-2x_1^{41}+3x_1^{40}-5x_1^{39}
 -7x_1^{38}+2x_1^{37}-3x_1^{36}+2x_1^{35}+2x_1^{34}-7x_1^{33}+6x_1^{32} \\
&\, +3x_1^{31}+3x_1^{29}+x_1^{28}-6x_1^{27}-3x_1^{26}+x_1^{25}+4x_1^{24}-7x_1^{23}-x_1^{22}+5x_1^{21}+3x_1^{20}
 -4x_1^{19}+2x_1^{18}-8x_1^{17}+5x_1^{16} \\
&\, +8x_1^{15}+2x_1^{14}+5x_1^{13}-4x_1^{12}+7x_1^{11}-2x_1^{10}-4x_1^{9}+4x_1^{8}
 +4x_1^{7}+3x_1^{6}-4x_1^{5}-5x_1^{4}-8x_1^{3}+x_1^2-x_1-1.
\end{align*}
}
The algebraic degree of this problem is $48$,
so the upper bound is sharp in this case.
\end{exm}

\subsection{$p$-th order cone programming}

The $p$-th order cone programming (pOCP) has the standard form
\be \label{pOCP}
\baray{rl}
\min_{x\in\re^n} &   c^Tx \\
s.t. &  a_i^Tx+b_i - \| C_ix+d_i \|_p \geq 0 , \, i=1,\cdots, \ell
\earay
\ee
where $c,a_i,b_i,C_i,d_i$ are matrices or vectors of appropriate dimensions.
This is also a convex optimization problem.
Let $x^*$ be one optimizer, and
assume the constraints with indices $1,\cdots,m$ are active at $x^*$.
Suppose the matrices $C_i$ has $r_i$ rows.
When some $r_i=1$, the constraint
$ a_i^Tx+b_i-\|C_ix+d_i\|_p \geq 0 $
is equivalent to two linear constraints
\[
-(a_i^Tx + b_i) \leq  C_i x + d_i \leq  a_i^Tx + b_i.
\]
Like the SOCP case, assume $1=r_{1}=\cdots=r_{k}<r_{k+1}\leq \cdots \leq r_{m}$.
Then problem~\reff{pOCP} is equivalent to
\begin{align*}
\min_{x\in\re^n} & \quad c^Tx \\
s.t. & \quad   a_i^Tx+b_i + \sig_i (C_ix+d_i) \geq 0 , \, i=1,\ldots, k \\
     & \quad  (a_i^Tx+b_i)^p - \sum_{j=1}^{r_i} (C_ix+d_i)_j^p  \geq  0 , \, i=k+1,\ldots, m
\end{align*}
where $\sig_i$ is chosen such that $a_i^Tx^* + b_i +\sig_i (C_i x^* + d_i)=0$.
In this situation
{\small
\[
D_{n-m}(0,\cdots,0,\underbrace{p-1,\cdots,p-1}_{m-k \mbox{ times } })
= \sum_{i_{k+1}+\cdots+i_m = n-m} (p-1)^{i_{k+1}+\cdots+i_m }
= (p-1)^{n-m} \binom{n-k-1}{m-k-1}.
\]
}
By Corollary~\ref{cor:deg},
the algebraic degree of $x^*$ is therefore bounded by
\be \label{eq:pdeg}
p^{m-k}(p-1)^{n-m} \binom{n-k-1}{m-k-1}.
\ee
When $k=m$, problem~\reff{pOCP} is reducible to some linear programming
and hence its algebraic degree is one.

Now we discuss the sharpness of degree bound~\reff{eq:pdeg}.
Similarly to the SOCP case,
for every $i=k+1,\cdots,m$, define the set of polynomials $U_{i}$ as
\[
U_i = \left\{
 (a_i^Tx+b_i)^p - \sum_{1\leq j \leq r_i}  \af_{j}^p (C_ix+d_i)_j^p :\,
\af_1,\cdots, \af_{r_i} \in \re
\right\}.
\]
Then define affine spaces $V_{i}$ as follows:
\begin{align*}
V_i &= \left\{
(a_i^Tx+b_i)^p - \sum_{1\leq j \leq r_i}  \bt_j (C_ix+d_i)_j^p :\,
\bt_1,\cdots, \bt_{r_i} \in \cpx
\right\}, \, i=k+1,\cdots, m.
\end{align*}
Then every set $U_{i}$ intersects any Zariski open subset of the affine space $V_{i}$.
On the other hand, the set of common zeros of the linear polynomials
$$a_i^Tx+b_i + \sig_i (C_ix+d_i),\, i=1,...,k$$
and all the polynomials in the union $\bigcup_{i=k+1}^m V_i$
is contained in the set $Z$ defined by \reff{eq:Zdef}.
Therefore, for generic choices of
$ a_i,b_i,C_i,d_i$, if $r_{k+1}+\cdots+r_m+m > n$,
the set $Z$ is empty,
and hence Remark~\ref{BAS} implies that the degree bound
given by formula~\reff{eq:pdeg} is sharp.

\begin{exm}
Consider the case $p=4$ and the polynomials
{\scriptsize
\begin{align*}
f_0 & = 9 x_1-5 x_2+3 x_3+2 x_4 \\
f_1 & = (1-6 x_1-6 x_2+4 x_3-9 x_4)^4
-(7-6 x_1+22 x_2-x_3+x_4)^4 -(11+x_1-x_2-8 x_3+3 x_4)^4 \\
& \quad -(-13+7 x_1+16 x_2-7 x_3+9 x_4)^4-(3-11 x_1+14 x_2-8 x_3+5 x_4)^4
-(8+9 x_1-10 x_2+2 x_3+2 x_4)^4.
\end{align*}
}
For the above polynomials, the inequality constraint must be active
since the objective is linear.
By formula~\reff{eq:pdeg}, the algebraic degree of the optimal solution is
bounded by $p^m (p-1)^{n-m} \binom{n-1}{m-1}=108$.
Symbolic computation shows
the optimal coordinate $x_1$ is a root of the univariate polynomial of degree $108$
(whose coefficients are modulo $17$)
{\tiny
\begin{align*}
&x_1^{108}-3 x_1^{107}-8 x_1^{106}+7 x_1^{105}+3 x_1^{104}-2 x_1^{103}-4 x_1^{102}-6 x_1^{101}+2 x_1^{100}
+8 x_1^{99}-8 x_1^{98}+5 x_1^{97}-3 x_1^{96}-3 x_1^{95} \\
&+4 x_1^{94}+3 x_1^{93}+7 x_1^{92}-4 x_1^{91}+6 x_1^{90} +x_1^{89}+7 x_1^{88}-x_1^{87}-5 x_1^{86}-6 x_1^{85}
+x_1^{84}+5 x_1^{83}-x_1^{81}+7 x_1^{80}+8 x_1^{79}\\
&-6 x_1^{78}+7 x_1^{77}+2 x_1^{76}-3 x_1^{75}+4 x_1^{74}-6 x_1^{73} -6 x_1^{72} +x_1^{70}+2 x_1^{69}-x_1^{68}
+8 x_1^{67}-3 x_1^{66}+5 x_1^{65}+4 x_1^{64}+x_1^{63}\\
&+x_1^{62}-2 x_1^{61}-x_1^{60}+3 x_1^{59}-7 x_1^{58}-7 x_1^{57}+7 x_1^{55}-3 x_1^{54}-3 x_1^{53}-8 x_1^{52}
-4 x_1^{51}-4 x_1^{50}-3 x_1^{49}-4 x_1^{48}+x_1^{47}\\
&+8 x_1^{46}+4 x_1^{45}-4 x_1^{44}-8 x_1^{43}-8 x_1^{42}-7 x_1^{41}-5 x_1^{40}+4 x_1^{39}-5 x_1^{38}-7 x_1^{37}
+4 x_1^{36} -2 x_1^{35}+x_1^{34}+6 x_1^{33}+6 x_1^{32}\\
&-7 x_1^{31}-3 x_1^{30}-5 x_1^{29}+7 x_1^{28}+3 x_1^{27}+6 x_1^{26}+2 x_1^{24}-8 x_1^{23}-8 x_1^{22}-4 x_1^{21}
+8 x_1^{20}+8 x_1^{19}-3 x_1^{18}+6 x_1^{17}-5 x_1^{16}\\
&-8 x_1^{15}+8 x_1^{14}+8 x_1^{13}+6 x_1^{12}-5 x_1^{10}+3 x_1^{9}+2 x_1^{8}-2 x_1^{7}+6 x_1^{6}+4 x_1^{5}
+7 x_1^{3}-8 x_1^{2}+4x_1+5.
\end{align*}
}
So the algebraic degree of this problem is $108$,
and the bound given by formula~\reff{eq:pdeg} is sharp.
\end{exm}

\bigskip
\noindent
{\bf Acknowledgement.}
The authors are grateful to Gabor Pataki, Richard Rimanyi and Bernd Sturmfels
for motivating this paper.

\end{document}